\documentclass[aos,MSNbibl,seceqn,citesort,dvips]{arximspdf}

\doi{10.1214/11-AOS899}
\volume{39}
\issue{4}
\pubyear{2011}
\firstpage{2074}
\lastpage{2102}

\makeatletter

\newtheorem{theorem}{Theorem}
\newtheorem{corollary}{Corollary}

\newproclaim{Example}{Example}

\makeatother

\begin{document}
\begin{frontmatter}

\title{Parametric or nonparametric? A parametricness index for model selection\thanksref{T1}}
\runtitle{A parametricness index}

\thankstext{T1}{Supported by NSF Grant DMS-07-06850.}

\begin{aug}
\author[A]{\fnms{Wei} \snm{Liu}\ead[label=e1]{william050@stat.umn.edu}} and
\author[A]{\fnms{Yuhong} \snm{Yang}\corref{}\ead[label=e2]{yyang@stat.umn.edu}\ead[label=u1,url]{http://www.stat.umn.edu/\textasciitilde yyang}}
\runauthor{W. Liu and Y. Yang}
\affiliation{University of Minnesota}
\address[A]{School of Statistics\\
University of Minnesota\\
313 Ford Hall \\
224 Church Street S.E. \\
Minneapolis, Minnesota 55455\\
USA\\
\printead{e1}\\
\hphantom{E-mail: }\printead*{e2}\\
\printead{u1}} 
\end{aug}

\received{\smonth{3} \syear{2010}}
\revised{\smonth{1} \syear{2011}}

%
\begin{abstract}
In model selection literature, two classes of criteria perform well
asymptotically in different situations: Bayesian information criterion
(BIC) (as a representative) is consistent in selection when the true
model is finite dimensional (parametric scenario); Akaike's information
criterion (AIC) performs well in an asymptotic efficiency when the true
model is infinite dimensional (nonparametric scenario). But there is
little work that addresses if it is possible and how to detect the
situation that a specific model selection problem is in. In this work,
we differentiate the two scenarios theoretically under some conditions.
We develop a measure, parametricness index (PI), to assess whether a
model selected by a potentially consistent procedure can be practically
treated as the true model, which also hints on AIC or BIC is better
suited for the data for the goal of estimating the regression function.
A~consequence is that by switching between AIC and BIC based on the PI,
the resulting regression estimator is simultaneously asymptotically
efficient for both parametric and nonparametric scenarios. In addition,
we systematically investigate the behaviors of PI in simulation and
real data and show its usefulness.\looseness=1
\end{abstract}

%

\end{frontmatter}

\section{Introduction}\label{intro}

When considering parametric models for data analysis, model selection
methods have been commonly used for various purposes. If one candidate
model describes the data really well (e.g., a physical law), it is
obviously desirable to identify it. Consistent model selection rules
such as BIC \cite{Schwarz78} are proposed for this purpose. In
contrast, when the candidate models are constructed to progressively
approximate an infinite-dimensional truth with a decreasing
approximation error, the main interest is usually on estimation and one
hopes that the selected model performs optimally in terms of a risk of
estimating a target function (e.g., the regression function). AIC
\cite{Akaike73} has been shown to be the right criterion from an
asymptotic efficiency and also a minimax-rate optimality views (see
\cite{Yang05} for references).

The question if we can statistically distinguish between parametric and
nonparametric scenarios motivated our
research. In this paper, for regression based on finite-dimensional
models, we develop a simple
parametricness index (PI) that has the following properties.

\begin{longlist}[(1)]
\item[(1)] With probability going to 1, PI separates typical parametric and
nonparametric scenarios.

\item[(2)] It advises on whether identifying the true or best candidate
model is feasible at the given sample size or not
by assessing if one of the models stands out as a stable parametric
description of the data.

\item[(3)] It informs us if interpretation and statistical inference based
on the selected model are
questionable due to model selection uncertainty.

\item[(4)] It tells us whether AIC is likely better than BIC for the data
for the purpose of estimating
the regression function.

\item[(5)] It can be used to approximately achieve the better estimation
performance of AIC and BIC for both
parametric and nonparametric scenarios.
\end{longlist}


In the rest of the \hyperref[intro]{Introduction}, we provide a
relevant background of model selection and present views on some
fundamental issues.

\subsection{Model selection criteria and their possibly conflicting properties}

To assess performance of model selection criteria, pointwise
asymptotic results (e.g.,
\cite
{Claeskens03,Fan01,IngWei05,Leeb08,Lv09,McQ98,Meinshausen06,Nishii84,Potscher91,Shao97,Shibata81,Stone79,Wang07,Yang07,Zhang08,Zhao06,Zou06})
have been established mostly in terms of either
selection consistency or an asymptotic optimality. It is well known
that AIC \cite{Akaike73}, $C_{p}$ \cite{Mallows73} and FPE
\cite{Akaike69,Shibata84} have an asymptotic
optimality property which says the accuracy of the estimator based on the
selected model is asymptotically the same as the best candidate model
when the
true model is infinite dimensional. In contrast,
BIC and the like are
consistent when the true model is finite dimensional and is among the
candidate models (see \cite{Shao97,Yang05} for references).

Another direction of model selection theory focuses on oracle risk bounds
(also called index of resolvability bounds). When the candidate models are
constructed to work well for target function classes, this approach yields
minimax-rate or near minimax-rate optimality results.
Publications of work in this direction include
\cite
{Baraud02,Barron94,Barron99,Barron91,Birge06,Bunea06,Candes07,Devroye97,Donoho96,Efr85,Koltchinskii06,Yang98},
to name a few. In
particular, AIC type of model selection methods are minimax-rate
optimal for both parametric and nonparametric scenarios under square
error loss for estimating
the regression function (see
\cite{Barron99,Yang05}). A remarkable feature of the works inspired by
\cite{Barron91}
is that with a complexity penalty (other
than one in terms of model dimension) added
to deal with a large number of (e.g., exponentially many) models,
the resulting risk or loss of the selected model automatically
achieves the best trade-off between approximation error, estimation
error and the model complexity, which provides tremendous theoretical
flexibility to
deal with a fixed countable list of models (e.g., for series expansion
based modeling) or a list of models chosen to depend on the
sample size (see, e.g., \cite{Barron99,Yang98,Yang99}).

While pointwise asymptotic results are certainly of interest, it
is not surprising that the limiting behaviors can be very different from
the finite-sample reality, especially when model selection is involved.
(see, e.g., \cite{Kabaila02,Danilov04,Leeb06}).

The general forms of AIC and BIC make it very clear that they and
similar criteria (such as GIC in \cite{RaoWu89}) cannot simultaneously
enjoy the properties of consistency in
a parametric scenario and asymptotic optimality in a nonparametric
scenario. Efforts
have been put on using penalties that are data-dependent and adaptive
(see,
e.g., \cite{BYY,George00,Hansen99,Ing07,ShenHuang06,Shen02,Yang072}).
Yang \cite{Yang072} showed that the asymptotic optimality of BIC for a
parametric scenario (which follows directly from consistency of BIC) and
asymptotic optimality of AIC for a nonparametric scenario can be shared
by an
adaptive model selection criterion. A similar two-stage adaptive model
selection rule for time series autoregression has been proposed by Ing
\cite{Ing07}. However, Yang \cite{Yang05,Yang072} proved that no
model selection procedure can be both consistent (or pointwise
adaptive) and minimax-rate optimal at
the same time.
As will be seen, if we can properly distinguish between parametric and
nonparametric
scenarios, a~consequent data-driven choice of AIC or BIC simultaneously
achieves asymptotic
efficiency for both parametric and nonparametric situations.

\subsection{Model selection: A gap between theory and practice}

It is well known
that for a typical regression problem with a number of predictors, AIC and
BIC tend to choose models of significantly different sizes, which may
have serious practical consequences.
Therefore, it is important to decide which criterion to apply for a
data set at
hand. Indeed, the conflict between AIC and BIC has received a lot of
attention not only in the statistics literature but also in fields such
as psychology and biology (see, e.g.,
\cite{Berger01,BA04,Chatfield95,Freedman95,Zhang97,Sober04}).
There has been a
lot of debate from not only
statistical but also philosophical perspectives, especially about the
existence of a true
model and the ultimate goal of statistical modeling.
Unfortunately, the current theories on model selection have little
to offer to address this issue. Consequently, it is rather common that
statisticians/statistical users resort to the ``faith'' that the true model
certainly cannot be finite dimensional for the choice of AIC, or to
the strong preference of parsimony or the goal of model identification
to defend his/her use of BIC.

To us, this disconnectedness between theory and practice of model selection
needs not to continue. From various angles, the question whether or not
AIC is more appropriate
than BIC for the data at hand should and can be addressed statistically rather
than based on one's preferred assumption. This is the major motivation
for us to try to go beyond
presenting a few theorems in this work.

We would like to quote a leading
statistician here:

``\textit{It does not seem helpful just to say that all models are
wrong. The very word model implies simplification and idealization. The idea
that complex physical, biological and sociological systems can be exactly
described by a few formulae is patently absurd. The construction of
idealized representations that capture important stable aspects of such
systems is, however, a vital part of general scientific analysis and
statistical models, especially substantive ones (Cox, 1990), do not seem
essentially different from other kinds of model}'' (Cox \cite{Cox95}).

Fisher in his pathbreaking 1922 paper \cite{Fisher22}, provided
thoughts on the foundations
of statistics, including model specification. He stated:
``More or less elaborate forms will be suitable according to the
volume of the data.'' Cook \cite{Cook07} discussed Fisher's insights in details.
%
%
%
%
%

We certainly agree with the statements by Fisher and Cox. What we are
interested in this and future work on model selection is to address the
general question that in
what ways and to what degrees a selected model is useful.

Finding a stable finite-dimensional model to describe the
nature of the data as well as to predict the future is very appealing.
Following up in the spirit of Cox mentioned above, if a model stably
stands out
among the competitors, whether it is the true model or not, from a
practical perspective,
why should not we extend the essence of consistency to mean
the ability to find it? In our view, if we are to accept any
statistical model (say infinite dimensional) as a useful vehicle to
analyze data, it is
difficult to philosophically reject the more restrictive assumption of
a finite-dimensional model, because both are convenient and
certainly simplified descriptions of the reality, their difference
being that between 50 paces and 100 paces as in the 2,000 year old
Chinese idiom
\textit{One who retreats fifty paces mocks one who retreats a hundred}.

The above considerations lead to the question: Can we construct
a practical measure that gives us a proper indication on whether the
selected model deserves to be crowned as the best model \textit{at the time
being}? We emphasize \textit{at the time being} to make it clear that we
are not going after the best limiting model (no matter how that is
defined), but instead we seek a
model that stands out for sample sizes around what we have now.

While there are many different performance measures that we can use to
assess if one model stands out, following our results on
distinguishing between parametric and nonparametric scenarios, we focus on
an estimation accuracy measure. We call it \textit{parametricness index} (PI), which is relative to the list of candidate
models and the sample size. Our theoretical results show that this
index converges to infinity for a parametric scenario and converges to
1 for a typical
nonparametric scenario. Our suggestion\vadjust{\goodbreak} is that when the index is
significantly larger than 1, we can treat the selected model as a
stably standing out model from the estimation perspective. Otherwise,
the selected model is just among a few or more equally well-performing
candidates. We call the former case practically parametric and the latter
practically nonparametric.

As will be demonstrated in our simulation work, PI can be close to 1
for a truly parametric scenario and large for a nonparametric scenario.
In our
view, this is not a problem. For instance, for a truly parametric scenario
with many small coefficients of various magnitudes, for a small or
moderate sample size, the
selected model will most likely be different from the true model and
it is also among multiple models that perform similarly in estimation
of the regression function. We would view this as ``practically nonparametric''
in the sense that with the information available we are not able to
find a single standing-out model and the model selected provides
a~good trade-off between approximation capability and model
dimension. In contrast, even if the true model is
infinite dimensional, at a given sample size, it is quite possible
that a number of terms are significant and others are too
small to be relevant at the given sample size. Then we are willing to
call it ``practically
parametric''
in the sense that as long as the sample size is not
substantially increased, the same model is expected to perform better
than the other candidates. For
example, in properly designed experimental studies, when a working model
clearly stands out and is very stable, then it is desirable to treat
it as a
parametric scenario even though we know surely it is an approximating
model. This is often the case in physical sciences when a law-like
relationship is evident under controlled experimental conditions. Note
that given an infinite-dimensional true
model and a list of candidate models, we may declare the selected
models to be practically parametric for some sample sizes and to be
practically nonparametric for others.

The rest of the paper is
organized as follows. In Section \ref{sec2}, we set up the regression framework and
give some notation. We then in Section \ref{sec3} develop the
measure PI and show that theoretically it differentiates a parametric scenario
from a nonparametric one under some conditions for both known and unknown
$\sigma^{2}$, respectively. Consequently, the pointwise asymptotic efficiency
properties of AIC and BIC can be combined for parametric and
nonparametric scenarios.
In Section~\ref{sec4}, we propose a proper use of PI for
applications. Simulation studies and real
data examples are reported in Sections \ref{sec5} and \ref{sec6}, respectively.
Concluding remarks are given in
Section \ref{sec7} and the proofs are in the \hyperref[app]{Appendix}.

\section{Setup of the regression problem}\label{sec2}

Consider the regression model
\[
Y_{i}=f(x_{i})+\varepsilon_{i},\qquad i=1,2,\ldots, n,
\]
where $x_{i}=(x_{i1},\ldots, x_{ip})$ is the value of a $p$-dimensional
fixed design variable at the $i$th observation, $Y_{i}$ is
the response, $f$ is the true regression\vadjust{\goodbreak} function, and the random
errors $%
\varepsilon_{i}$ are assumed to be independent and normally distributed with
mean zero and variance $\sigma^{2}>0$.

To estimate the regression function, a list of linear models are being
considered, from which one is to be selected:
\[
Y=f_{k}(x,\theta_{k})+\varepsilon^{\prime},
\]
where, for each $k$, $\mathcal{F}_{k}=\{f_{k}(x,\theta_{k}), \theta_{k}
\in
\Theta_{k}\}$ is a family of regression functions linear in the parameter
$\theta_{k}$ of finite dimension $m_{k}$.
Let $\Gamma$ be the collection of the model
indices $k$. $\Gamma$ can be fixed or change with the sample size.

The above framework includes the usual subset-selection and order-selec\-tion
problems in linear regression. It also includes nonparametric regression
based on series expansion, where the true function is approximated by linear
combinations of appropriate basis functions, such as polynomials,
splines or
wavelets.

Parametric modeling typically intends to capture the essence of the
data by
a finite-dimensional model, and nonparametric modeling tries to achieve the
best trade-off between approximation error and estimation error for a~target
infinite-dimensional function. See, for example,
\cite{Yang992} for the general relationship between rate of convergence
for function estimation and full or sparse approximation based on a linear
approximating system.

Theoretically speaking, the essential difference between parametric and
nonparametric scenarios in our context is that the best model has no
approximation
error for the former and all the candidate models have nonzero
approximation errors for the latter.

In this paper, we consider the least squares estimators when defining
the parametricness index, although
the model being examined can be based any consistent model selection
method that may
or may not involve least squares estimation.

\subsection*{Notation and definitions}

Let ${\mathbf Y}_{n}=(Y_{1},\ldots, Y_{n})^{T}$ be the response vector and
$M_{k}$ be
the projection matrix for model $k$. Denote $\hat{{\mathbf Y}}_{k}=M_{k}{\mathbf Y}_{n}$.
Let $f_{n}=(f(x_{1}), \ldots,
f(x_{n}))^{T}$, $%
e_{n}=(\varepsilon_{1}, \ldots, \varepsilon_{n})^{T}$, and $I_{n}$ be
the identity
matrix. Let \mbox{$\|\cdot\|$} denote the Euclidean
distance in the $R^{n}$ space, and let $%
\operatorname{TSE}(k)=\|f_{n}-\hat{{\mathbf Y}}_{k}\|^{2}$ be
the total square error of the LS estimator from model $k$.

Let the
rank of $M_{k}$ be $r_{k}$. In this work, we do not assume that all the
candidate models have the
rank of the design matrix equal the model dimension~$m_k$, which may
not hold when a large number of models are considered. Let~$N_{j}$
denote the number of models with $r_{k}=j$ for $k\in\Gamma$. For a~given
model~$k$, let~$S_{1}(k)$ be the set of all sub-models $k^{\prime}$
of~$k$ in~$\Gamma$ such that $r_{k^{\prime}}=r_{k}-1$. Throughout the
paper, for technical convenience, we assume~$S_{1}(k)$ is not empty for
all $k$ with $r_{k}>1 $.

For a sequence $\lambda_{n} \geq(\log n)^{-1}$ and a constant $d\ge
0$, let
\[
IC_{\lambda_{n}, d }(k)=\|{\mathbf Y}_{n}-\hat{{\mathbf Y}}_{k}\|^{2}+\lambda
_{n}\log(n)r_{k}\sigma^{2}-n\sigma^{2}+d n^{1/2}\log(n)\sigma^{2},
\]
when $\sigma$ is known, and
\[
IC_{\lambda_{n}, d}(k,\hat{\sigma}^{2})=\|{\mathbf Y}_{n}-\hat{{\mathbf Y}}_{k}\|^{2}+\lambda_{n}\log(n)r_{k}\hat{\sigma}^{2}-n\hat{\sigma
}^{2}+d n^{1/2}\log(n)\hat{\sigma}^{2},
\]
when $\sigma$ is estimated by $\hat{
\sigma}$. A discussion on choice of $\lambda_{n}$ and $d$ will be given
later in Section \ref{sec35}. We emphasize that our use of $IC_{\lambda
_{n}, d}(k)$ or $IC_{\lambda_{n}, d}(k,\hat{\sigma}^{2})$ is for
defining the parametricness index as below and it may not be the one
used for model selection.

\section{Main theorems}\label{sec3}

Consider a potentially consistent model selection me\-thod (i.e., it
will select the true model with
probability going to 1 as $n\rightarrow\infty$ if the true model is
among the candidates). Let $\hat{k}_{n}$ be the
selected model at sample size $n$. We define the \textit{parametricness
index} (PI) as follows:
\begin{longlist}[(1)]
\item[(1)] When $\sigma$ is known,
\[
\mathrm{PI}_{n}= \cases{
\displaystyle \inf_{k\in S_{1}(\hat{k}_{n} )} \frac{IC_{\lambda_{n}, d
}(k)}{IC_{\lambda_{n}, d}(\hat{k}_{n})}, &\quad if
$r_{\hat{k}_{n}}>1$, \vspace*{2pt}\cr
n, &\quad if $r_{\hat{k}_{n}}=1$.}
\]
\item[(2)] When $\sigma$ is estimated by
$\hat{\sigma}$,
\[
\mathrm{PI}_{n}= \cases{
\displaystyle\inf_{k\in S_{1}(\hat{k}_{n} )}
\frac{IC_{\lambda_{n}, d}(k,\hat{\sigma}^{2})}{IC_{\lambda_{n},
d}(\hat{k}_{n},\hat{\sigma}^{2})}, &\quad if $r_{\hat
{k}_{n}}>1$, \vspace*{2pt}\cr
n, &\quad if $r_{\hat{k}_{n}}=1$.}
\]
\end{longlist}
The reason behind the definition is that a correctly specified
parametric model must be
very different from any sub-model (bias of a sub-model is dominatingly
large asymptotically speaking), but for a nonparametric scenario, the
model selected is only slightly affected in terms of estimation
accuracy when one or a few
least important terms are dropped. When $r_{\hat{k}_{n}}=1$, the value
of PI is arbitrarily defined as long as
it goes to infinity as $n$ increases.


\subsection{Parametric scenarios}\label{sec31}

Now consider a \textit{parametric scenario}: the true model at sample size
$n$ is in $\Gamma$ and denoted
by $k_{n}^{*}$ with $r_{k_{n}^{*}}$ assumed to be larger than 1. Let
$A_{n}=\inf_{k\in S_{1}(k_{n}^{*} )}
\|(I_{n}-M_{k})f_{n}\|^{2}/\sigma^{2}. $
Note that $A_{n}/n$ is the best approximation error (squared bias) of
models in $S_{1}(k_{n}^{*} )$.

\textit{Conditions}:
\begin{longlist}[(P1)]
\item[(P1)] There exists $0<\tau\le\frac{1}{2}$ such that $A_{n}$ is of order
$n^{{1/2}+\tau}$ or higher.
%
\item[(P2)] The dimension of the true model does not grow too fast with
sample size $n$ in the sense that $ r_{k_{n}^{*}}\lambda_{n}\log
(n)=o(n^{{1/2}+\tau})$.
\item[(P3)] The selection procedure is consistent: $P(\hat{k}_{n}=k_{n}^{*})
\rightarrow1$ as $n \rightarrow\infty$.
\end{longlist}
\begin{theorem} \label{thm1}
Assume conditions \textup{(P1)--(P3)} are satisfied for the parametric scenario.

\begin{longlist}
\item
With $\sigma^{2}$ known, we have
\[
\mathrm{PI}_{n} \stackrel{p}{\longrightarrow} \infty
\qquad\mbox{as $n \rightarrow\infty$}.
\]

\item When $\sigma$ is unknown, let $\hat{\sigma}^{2}_{n}=\frac{\|{\mathbf Y}_{n}-\hat{{\mathbf Y}}_{\hat{k}_n}\|^{2}}{n-r_{\hat{k}_n}}$. We also
have
\[
\mathrm{PI}_{n} \stackrel{p}{\longrightarrow} \infty
\qquad\mbox{as $n \rightarrow\infty$}.
\]
\end{longlist}
\end{theorem}

\textit{Remarks}:
(1) The
conditions (P1) basically eliminates the case that the true model and a
sub-model with one fewer term are not distinguishable with the
information available in the sample.

(2) In our formulation, we considered comparison of two immediately
nested models. One can consider comparing two nested
models with size difference $m$ ($m >1$) and similar results hold.

(3) The case $\lambda_{n}=1$ corresponds to using BIC in defining the
PI. And $\lambda_{n}=2/\log(n)$ corresponds to using AIC.

\subsection{Nonparametric scenarios}\label{sec32}

Now the true model at each sample size~$n$ is not in the list $\Gamma$
and may change with sample size, which we call a~\textit{nonparametric scenario}.
For $j<n$, denote
\[
B_{j,n}=\inf_{k\in\Gamma} \bigl\{ \bigl(\lambda_{n}\log(n)-1\bigr)j+\|(I_{n}-M_{k
})f_{n}\|^{2}/\sigma^{2} + d n^{1/2}\log(n)\dvtx r_{k }=j\bigr\},
\]
%
where the infimum is taken over all the candidate models with $r_{k}=
j$. For $1<j<n$, let $L_{j}=\max_{k \in\Gamma}\{ \operatorname{card}(S_{1}(k))\dvtx
r_{k}=j \}$. Let $P_{k^{(s)},k}= M_{k}-M_{k^{(s)}} $ be the difference
between the projection matrices of the two nested models. Clearly,
$P_{k^{(s)},k}$ is the projection matrix onto the orthogonal complement
of the column space of model $k^{(s)}$ with respect to that of the
larger model $k$.

\textit{Conditions}: There exist two sequences of integers $1\le a_{n}<
b_{n}<n$ (not necessarily known)
with $a_{n}\rightarrow\infty$ such that the following holds:
\begin{longlist}[(N1)]
\item[(N1)] $P(a_{n}\le r_{\hat{k}_{n}}\le b_{n})\rightarrow1$ and $\sup
_{a_{n}\le j \le b_{n}} \frac{B_{j,n}}{n-j} \rightarrow0 $ as $n
\rightarrow\infty$.


\item[(N2)]
There exist a positive sequence $ \zeta_{n} \rightarrow0$ and
constants $c_{0}>0$ such that for $a_{n}\le j \le b_{n}$,
\[
N_{j}\cdot L_{j}\le c_{0}e^{\zeta_{n} B_{j,n}},\qquad N_{j} \le
c_{0}e^{{B_{j,n}^{2}}/({10(n-j)})}
\]
and
\[
\limsup
_{n\rightarrow\infty}\sum_{j=a_{n}}^{b_{n}} e^{-
{B_{j,n}^{2}}/({10(n-j)})} = 0.
\]

\item[(N3)]
\begin{eqnarray*}\nonumber
\!\limsup_{n\rightarrow\infty} [\sup_{\{k\dvtx a_{n}\le r_{k}\le
b_{n}\}}\!\frac{\inf_{k^{(s)}\in S_{1}(k)}\|P_{k^{(s)},k}f_{n}\|^{2}}
{(\lambda_{n}\log(n)\,{-}\,1)r_{k}\,{+}\,\|(I_{n}\,{-}\,M_{k})f_{n}\|^{2}/\sigma^{2}\,{+}\,d n^{1/2}\log(n)}]
\,{=}\,0.\nonumber
\end{eqnarray*}
\end{longlist}
\begin{theorem} \label{thm2}
Assuming conditions \textup{(N1)--(N3)} are satisfied for a nonparametric scenario
and $\sigma^{2}$ is known, then we have
\[
\mathrm{PI}_{n} \stackrel{p}{\longrightarrow} 1
\qquad\mbox{as $n \rightarrow\infty$}.
\]
\end{theorem}

\textit{Remarks}:
(1) For nonparametric regression, for
familiar model selection methods, the order of $ r_{\hat{k}_{n}}$
can be\vspace*{1pt} identified (e.g., \cite{Ing07,Yang992}), sometimes loosing a
logarithmic factor, and (N1) is satisfied in a typical nonparametric situation.

(2) Condition (N2) basically ensures that the number of subset
models of each dimension does not grow too fast relative to
$B_{j,n}$. When the best model has a slower rate of convergence in
estimating $f$,
more candidate models can be allowed without detrimental selection
bias.

(3) Roughly speaking, condition (N3) says that when the model dimension
is in a range
that contains the selected model with probability approaching~1,
the least significant term in the regression function projection is
negligible compared to the sum of approximation error,\vspace*{1pt} the dimension of
the model times $\lambda_{n}\log(n)$, and the term $d n^{1/2}\log(n)$.
This condition is mild.

(4) A choice of $d>0$ can handle situations where the approximation
error decays fast, for example, exponentially fast (see Section \ref{sec34}), in
which case the stochastic fluctuation of $IC_{\lambda_{n}, d }$
with $d=0$ is relatively too large for PI to converge to 1 in
probability. In applications, for separating reasonably distinct
parametric and nonparametric scenarios, we recommend the choice of
$d=0$.

When $\sigma^{2}$ is unknown but estimated from the selected model,
$\mathrm{PI}_{n}$ is correspondingly defined. For $j<n$, let $E_{j,n}$ denote
\[
\inf_{k\in\Gamma, r_{k}=j } \bigl\{ \bigl[\bigl(\lambda_{n}\log(n)-1\bigr)j+d
n^{1/2}\log(n)\bigr] \bigl[1+\|(I_{n}-M_{k})f_{n}\|^{2}/\bigl((n-j)\sigma
^{2}\bigr) \bigr] \bigr\}.
\]


\textit{Conditions}: There exist two sequences of integers $1\le a_{n}< b_{n}<n$
with $a_{n}\rightarrow\infty$ such that the following holds.
\begin{longlist}[(N2$'$)]

%
\item[(N2$'$)]
There exist a positive sequence $\rho_{n}\rightarrow0$ and a constant
$c_{0}>0$ such that for $a_{n}\le j \le b_{n}$, $ N_{j}\cdot L_{j}\le
c_{0}e^{\rho_{n} E_{j,n}}$, and $\limsup_{n\rightarrow\infty}\sum
_{j=a_{n}}^{b_{n}} e^{- \rho_{n}E_{j,n} } = 0$.\vspace*{1pt}

\item[(N3$'$)] $ \limsup_{n\rightarrow\infty}[\sup_{\{
k\dvtx a_{n}\le r_{k}\le b_{n} \} }(({\inf_{k^{(s)}}\|
P_{k^{(s)},k}f_{n}\|^{2}})([(\lambda_{n}\log(n)-\break1)r_{k}+ d n^{1/2}\log
(n)][1+\|(I_{n}-M_{k})f_{n}\|^{2}/(\sigma^{2}(n-r_{k}))] )^{-1})] =0 $.
\end{longlist}
\begin{theorem} \label{thm3}
Assuming conditions \textup{(N1)}, \textup{(N2$'$)} and \textup{(N3$'$)}
hold for a~nonparametric scenario, then we have
\[
\mathrm{PI}_{n} \stackrel{p}{\longrightarrow} 1 \qquad\mbox{as $n
\rightarrow\infty$}.
\]
\end{theorem}

\subsection{PI separates parametric and nonparametric scenarios}

The results in Sections \ref{sec31} and \ref{sec32} imply that starting with a
potentially consistent model
selection procedure (i.e., it will be consistent if one of the
candidate models holds),
the $\mathrm{PI}$ goes to $\infty$ and $1$ in probability in parametric
and nonparametric scenarios, respectively.
\begin{corollary}
Consider a model selection setting where $\Gamma_{n}$ includes models
of sizes approaching $\infty$ as $n \rightarrow\infty$. Assume the
true model is either parametric or nonparametric satisfying \textup{(P1)} and
\textup{(P2)}
or \textup{(N1)--(N3)}, respectively. Then $\mathrm{PI}_{n}$ has distinct limits in
probability for the two scenarios.
%
\end{corollary}


\subsection{Examples}\label{sec34}

We now take a closer look at the conditions (P1)--(P3) and (N1)--(N3)
for two settings: all subset selection and order selection (i.e., the
candidate models are nested).

(1) \textit{All subset selection}.

Let $p_{n}$ be the number of terms to be considered.

\begin{longlist}
\item
Parametric with true model $k_{n}^{*}$ fixed.

In this case, $A_{n}$ is typically of order $n$ for a reasonable design
and then condition (P1) is met. Condition (P2) is obviously satisfied
when $\lambda_{n}=o(n^{1/2})$.

\item Parametric with $k_{n}^{*}$ changing: $r_{k_{n}^{*}}$ increases
with $n$.

In this case, both $r_{k_{n}^{*}}$ and $p_{n}$ go to infinity with $n$.
Since there are more and more terms in the true model, in order for
$A_{n}$ not to be too small, the terms should not be too highly
correlated. An extreme case is that one term in the true model is
almost linearly dependent on the others. Then $A_{n} \thickapprox0$.
To understand condition~(P1) in terms of the coefficients in the true
model, under an orthonormal design, condition~(P1) is more or less
equivalent to  the square of the smallest coefficient in the true
model is of order $n^{\tau-1/2}$ or higher. Since $ \tau$ can be
arbitrarily close to 0, the smallest coefficient should basically be
larger than $n^{-{1/4}}$.

\item Nonparametric.

Condition (N1) holds for any model selection method that yields a
consistent regression estimator of $f$.
The condition $N_{j} \le c_{0}e^{{ B_{j,n}^{2}}/({10(n-j)})} $ is
roughly\vspace*{1pt} equivalent to
$j\log(p_{n}/j)\le[dn^{1/2}\log(n)+\lambda_{n}\log(n)j+\|
(I_{n}-M_{k})f_{n}\|^{2}/\sigma^{2}]^{2}/\break10(n-j) $ for $a_{n}\le j \le
b_{n}$. A sufficient condition is $p_{n}\le b_{n}e^{
B_{j,n}^{2}/(10(n-j)b_{n})}$ for $a_{n}\le j \le b_{n}$. As to the
condition $N_{j} \cdot L_{j}\le c_{0}e^{ \zeta{'}_{n}B_{j,n} } $, as
long as\break $\sup_{a_{n}\le j \le b_{n}}\frac{B_{j,n}}{n-j}\rightarrow0$,
then it is implied by the above one. For the condition $\sum
_{j=a_{n}}^{b_{n}} e^{-{ B_{j,n}^{2}}/({10(n-j)})}\rightarrow0 $, it
is automatically satisfied for any $d>0$ and also satisfied for $d=0$
when the approximation error does not decay too fast.\looseness=-1
\end{longlist}

%
%
%


(2) \textit{Order selection in series expansion}.

%
%
%

We only need to discuss the nonparametric scenario. (The parametric
scenarios are similar to the above.)

In this setting, there is only one model of each dimension. So
condition~(N2) reduces to: $\sum_{j=a_{n}}^{b_{n}} e^{-{
B_{j,n}^{2}}/({10(n-j)})}\rightarrow0 $. Note that\break $\sum
_{j=a_{n}}^{b_{n}} e^{-{ B_{j,n}^{2}}/({10(n-j)})}<
(b_{n}-a_{n})\cdot e^{- (\log(n))^{2}/10}< n\cdot e^{- (\log
(n))^{2}/10} \rightarrow0$.

To check condition (N3), for a demonstration, consider orthogonal
designs. Let ${\bolds\Phi} =\{\phi_{1}(x),\ldots, \phi_{k}(x), \ldots\}
$ be a collection of orthonormal basis functions and the true
regression function is
$f(x)=\sum_{i=1}^{\infty} \beta_{i}\phi_{i}(x)$. For model $k$, the
model with the first $k$ terms, $\inf_{k^{(s)}\in S_{1}(k)} \|
P_{k^{(s)},k}f_{n}\|^{2}$ is roughly $\beta_{k}^{2}\|\phi_{k}({\mathbf X})\|
^{2}$ and $\|(I_{n}-M_{k})f_{n}\|^{2}$ is roughly $ \sum_{i=k+1}^{\infty
} \beta_{i}^{2}\|\phi_{i}({\mathbf X})\|^{2} $, where $\phi_{i}({\mathbf
X})=(\phi_{i}(x_{1}),\break
\ldots, \phi_{i}(x_{n}))^{T}$. Since $\|\phi
_{i}({\mathbf X})\|^{2}$ is of order $n$, condition (N3) is roughly
equivalent to the following:
\[
\limsup_{n\rightarrow\infty}\biggl[ \sup_{a_{n}\le k \le b_{n}} \frac
{n\beta_{k}^{2}}{(\lambda_{n}\log(n)-1)k+n\sum_{i=k+1}^{\infty} \beta
_{i}^{2}/\sigma^{2}+dn^{1/2}\log(n)} \biggr] = 0.
\]
Then a sufficient condition for condition (N3) is that $d=0$ and
\[
\lim
_{k\rightarrow\infty}\frac{\beta_{k}^{2}}{\sum_{i=k+1}^{\infty} \beta
_{i}^{2} }=0,
\]
which is true if $\beta_{k} = k^{-\delta} $ for some
$\delta>0$ but not true if $\beta_{k} = e^{-ck} $ for some $c>0$.
When $\beta_{k}$ decays faster so that $\frac{\beta_{k}^{2}}{\sum
_{i=k+1}^{\infty} \beta_{i}^{2} }$ is bounded away from zero and
$\sup_{a_{n}\le k \le b_{n}}|\beta_{k}|=o(\frac{\sqrt{\log
(n)}}{n^{1/4}} ) $, any choice of $d>0$
makes condition (N3) satisfied. An example is the exponential-decay
case, that is, $\beta_{k} = e^{-ck} $ for some $c>0$. According to \cite
{Ing07}, when $\hat{k}_{n}$ is selected by BIC for order selection, we
have that $r_{\hat{k}_{n}}$ basically falls within a constant from
$\frac{1}{2c}\log(n/\log(n))$ in probability. In this case, $\beta
_{k}\approx\frac{\sqrt{\log(n)}}{n^{1/2}} $ for $k\approx\frac
{1}{2c}\log(n/\log(n)) $. Thus, condition (N3) is satisfied.

%
%
%

\subsection{\texorpdfstring{On the choice of $\lambda_n$ and $d$}{On the choice of lambda n and d}} \label{sec35}

A natural choice of $(\lambda_{n},d)$ is $\lambda_{n}=1$ and
$d=0$, which is expected to work well to distinguish parametric and
nonparametric scenarios that are not too close to each other for order
selection or
all subset selection with $p_n$ increasing not fast in $n$.
Other choices can handle more difficult situations, mostly
entailing the satisfaction of (N2) and~(N3). With a
larger $\lambda_{n} $ or~$d$, $\mathrm{PI}$ tends to be closer to 1 for
a nonparametric case, but at the same time, it makes a parametric case
less obvious. When there are many models being considered, $\lambda_{n}$
should not be too small so as to avoid severe selection bias. The
choice of $d>0$ handles fast decay of the approximation error in
nonparametric scenarios, as mentioned already.

%

\subsection{Combining strengths of AIC and BIC}

From above, for any given cutoff point bigger than 1, the $\mathrm{PI}$ in a
parametric scenario will eventually exceed it while the $\mathrm{PI}$ in a
nonparametric scenario will eventually drops below it when the sample
size gets large enough.

It is well known that AIC is asymptotically loss (or risk) efficient
for nonparametric scenarios and
BIC is consistent when there are fixed finite-dimensional correct
models, which implies that BIC is asymptotically loss efficient \cite{Shao97}.
\begin{corollary}
For a given number $c>1$, let $\delta$ be the model selection procedure
that chooses either
the model selected by AIC or BIC as follows:
\[
\delta=
\cases{
\mbox{AIC}, &\quad if $\mathrm{PI}<c$, \cr
\mbox{BIC}, &\quad if $\mathrm{PI}\ge c$.}
\]
Under conditions \textup{(P1)--(P3)/(N1)--(N3)}, $\delta$ is asymptotically loss
efficient in both parametric and nonparametric scenarios
as long as AIC and BIC are loss efficient for the respective scenarios.
\end{corollary}

\textit{Remarks}: (1) Previous work on sharing the strengths of
AIC and BIC utilized
minimum description length criterion in an adaptive fashion
\cite{BYY,Hansen99}, or
flexible priors in a Bayesian framework \cite{George00,Erven08}.
Ing \cite{Ing07} and Yang \cite{Yang072} established (independently)
simultaneous
asymptotic efficiency for both parametric and nonparametric
scenarios.

(2) Recently, Erven, Gr\"{u}nwald and de Rooij \cite{Erven08} found that
if a cumulative risk (i.e., the sum of risks from the sample size $1$
to $n$)
is considered instead of the usual risk at sample size $n$, then the
conflict between consistency in selection and minimax-rate optimality
shown in
\cite{Yang05}
can be resolved by a~Baye\-sian strategy that allows switching between models.

\section{PI as a model selection diagnostic measure,
that is, Practical Identifiability of the best model}\label{sec4}

Based on the theory presented in the previous section, it is natural to
use the simple rule for
answering the question if we are in a parametric or nonparametric
scenario: call it parametric if
PI is larger than $c$ for some $c >1$ and otherwise nonparametric.
Theoretically speaking, we will be right with
probability going to one.

Keeping in mind that the concepts such as parametric, nonparametric,
consistency and asymptotic efficiency
are all mathematical abstractions that hopefully characterize the
nature of the data and the behaviors
of estimators at the given sample size, our intended use of PI is not a
rigid one so as to be practically
relevant and informative, as we explain below.

Both parametric and nonparametric methods have been widely used in
statistical applications. One specific approach to nonparametric
estimation is to use parametric models as approximations to an
infinite-dimensional function, which is backed up by
approximation theories. However, it is in this case that the boundary
between parametric and nonparametric
estimations becomes blurred, and our work tries to address the
issue.

From a theoretical perspective, the difference between parametric and
nonparametric modeling is quite clear in this context. Indeed, when one
is willing to
assume that the data come from a member in a parametric family, the
focus is then naturally on the estimation of the parameters, and
finite-sample and large sample properties (such as UMVUE, BLUE,
minimax, Bayes and asymptotic efficiency) are well understood. For
nonparametric estimation, given infinite-dimensional smooth function
classes, various approximation systems (such as polynomial,
trigonometric and wavelets) have been shown to lead to minimax-rate
optimal estimators via various statistical methods (e.g.,
\cite{Birge86,Donoho96,IH77,Stone82}). In addition, given a function
class defined in terms of approximation error decay behavior by an
approximating system, rates of convergence of minimax risks have been
established (see, e.g., \cite{Yang992}). As is expected, the
optimal model size (in rate) based on linear approximation depends on
the sample
size (and other things) for a nonparametric scenario. In particular, for
full and sparse approximation sets of functions, the minimax theory
shows that for a typical nonparametric scenario, the optimal model
size makes the approximation error (squared bias) roughly equal to
estimation
error (model dimension over the sample size) \cite{Yang992}.
Furthermore, adaptive estimators that are simultaneously optimal for
multiple function classes can be obtained by model selection or
model combining (see, e.g., \cite{Barron99,Yang00} for many references).

From a practical perspective, unfortunately, things are much less clear.
Consider, for example, the simple case of polynomial regression. In
linear regression textbooks, one often finds data that show obvious
linear or quadratic behavior, in which case perhaps most statisticians
would be
unequivocally happy with a linear or quadratic model (think of Hooke's
law for describing elasticity). When the underlying
regression function is much more complicated so as to require 4th or
5th power, it becomes difficult to classify the situation as
parametric or nonparametric. While few (if any) statisticians would
challenge the notion that in both cases, the model is only an
approximation to reality, what makes the difference in calling one case
parametric quite comfortably but not the other?
Perhaps simplicity and stability of the model play key roles as mentioned
in Cox \cite{Cox95}. Roughly
speaking, when
a model is simple and fits the data excellently (e.g., with $R^2$ close
to 1)
so that there is little room to significantly improve the fit, the model
obviously stands out.
In contrast, if we have to use a 10th order polynomial to be able to fit
the data with 100 observations, perhaps
few would call it a parametric scenario. Most of the situations may
be in between.

Differently from the order selection
problem, the case of subset selection in regression is substantially
more complicated
due to the much increased complexity of the list of models. It seems to
us that
when all subset regression is performed, it is usually
automatically treated as
a parametric
problem in the literature. While this is not surprising, our view is
different. When the
number of variables
is not very small relative to the sample size and the error variance, the
issue of model selection does not seem to be too different from order
selection for polynomial regression where a high polynomial power is
needed. In our view, when analyzing data (in contrast to asymptotic
analysis), if one explores over a number of parametric models, it is not
necessarily proper to treat the situation as a parametric one (i.e.,
report standard errors
and confidence intervals for parameters and make interpretations based on
the selected model without assessing its reliability).

Closely related to the above discussion is the issue of model
selection uncertainty (see, e.g., \cite{Breiman96,Chatfield95}).
It is an important issue to know when we
are in a situation where a relatively simple and reliable model stands
out in a proper sense and thus can be used as the ``true'' model for
practical purposes, and when a selected model is just one out of
multiple or even many possibilities among the candidates
at the given sample size. In the
first case, we would be willing to call it parametric (or more
formally, practically parametric) and the latter (practically)
nonparametric.

We should emphasize that in our review, our goal is not exactly finding out
whether the underlying model is finite dimensional (relative to the
list of
candidate models) or not. Indeed, we will not be unhappy to declare a
truly parametric scenario nonparametric when around the current sample size
no model selection criterion can possibly identify it with confidence
and then take
advantage of it, in which case, it seems better to view the models as
approximations to the true one and we are just making a tradeoff
between the approximation error and estimation error. In contrast, we
will not be shy to continue calling a truly nonparametric model
parametric should we be given that knowledge by an oracle if one model
stands out at the current sample size and the contribution of the
ignored features is so small that it is clearly better to be ignored at
the time being. When the sample size is much increased, the enhanced
information allows discovery of the relevance of some additional
features and then we may be in a practical nonparametric scenario. As the
sample size further increases, it may well be that a parametric model
stands out until reaching a larger sample size where we enter a
practical nonparametric scenario again, and so on.

Based on hypothesis testing theories, obviously, at a given sample size,
for any true parametric distribution in one of the candidate families
from which the data are generated, one has a nonparametric
distribution (i.e., not in any of the candidate families) that cannot be
distinguished from the true distribution. From this perspective,
pursuing a rigid finite-sample distinction between parametric
and nonparametric scenarios is improper.

PI is relative to the list of candidate models and the sample
size. So it is perfectly possible (and fine) that for one list of
models, we declare the situation to be parametric, but for a different
choice of candidate list, we declare nonparametriness.

\section{Simulation results}\label{sec5}

In this section, we consider single-predictor and multiple-predictor
cases, aiming at a serious understanding of the practical utility of
PI. In all the numerical examples in this paper, we choose $\lambda
_{n}=1$ and $d=0$.

\subsection{Single predictor}

\begin{Example}\label{Example1}
Compare two different situations:

Case 1: $Y=3\sin(2\pi x)+\sigma_{1} \varepsilon$.

Case 2: $ Y=3-5x+2x^{2}+1.5x^{3}+0.8x^{4}+\sigma_{2} \varepsilon$,
where $\varepsilon\sim N(0, 1)$ and $x \sim N(0, 1)$.

BIC is used to select the order of polynomial regression between 1 and~30.
The estimated $\sigma$ from the selected model is used to calculate
the PI.


Quantiles for the PIs in both scenarios based on 300 replications are
%
%
\begin{table}[b]
\caption{Percentiles of PI for Example \protect\ref{Example1}}
\label{table1}
\begin{tabular*}{\tablewidth}{@{\extracolsep{\fill}}lcccccc@{}}
\hline
& \multicolumn{3}{c}{\textbf{Case 1}} &
\multicolumn{3}{c@{}}{\textbf{Case 2}} \\[-4pt]
& \multicolumn{3}{c}{\hrulefill} & \multicolumn{3}{c@{}}{
\hrulefill}\\
\textbf{Percentile} & \textbf{Order selected} & \textbf{PI}
& $\bolds{\hat{\sigma}}$ & \textbf{Order selected} &
\textbf{PI} & $\bolds{\hat{\sigma}}$ \\
\hline
10\% & \hphantom{0}1 & 0.47 & 2.78 & 4 & 1.14 & 6.53 \\
20\% & 13 & 1.02 & 2.89 & 4 & 1.35 & 6.67 \\
50\% & 15 & 1.12 & 3.03 & 4 & 1.89 & 6.96 \\
80\% & 16 & 1.34 & 3.21 & 4 & 3.15 & 7.31 \\
90\% & 17 & 1.54 & 3.52 & 4 & 4.21 & 7.49 \\
\hline
\end{tabular*}
\end{table}
presented in Table \ref{table1}.
\end{Example}
\begin{Example}\label{Example2}
Compare the following two situations:\vspace*{1pt}

Case 1: $Y=1-2x+1.6x^{2}+0.5x^{3}+3\sin(2\pi x)+\sigma\varepsilon$.

Case 2: $ Y=1-2x+1.6x^{2}+0.5x^{3}+\sin(2\pi x)+\sigma\varepsilon$.

The two mean functions are the same except the coefficient of the
$\sin(2\pi x)$ term.
As we can see from Table \ref{table2}, although both cases are of a
nonparametric nature, they have different behaviors in terms of model
selection uncertainty and PI values. Case 2 can be called ``practically''
parametric and the large PI values provide information in this regard.

%
\begin{table}
\caption{Percentiles of PI for Example \protect\ref{Example2}}
\label{table2}
\begin{tabular*}{\tablewidth}{@{\extracolsep{\fill}}lcccccc@{}}
\hline
& \multicolumn{3}{c}{\textbf{Case 1}} &
\multicolumn{3}{c@{}}{\textbf{Case 2}} \\[-4pt]
& \multicolumn{3}{c}{\hrulefill} & \multicolumn{3}{c@{}}{
\hrulefill}\\
\textbf{Percentile} & \textbf{Order selected} & \textbf{PI}
& $\bolds{\hat{\sigma}}$ & \textbf{Order selected} &
\textbf{PI} & $\bolds{\hat{\sigma}}$ \\
\hline
10\% & 15 & 1.01 & 1.87 & 3 & 1.75 & 1.99 \\
20\% & 15 & 1.05 & 1.92 & 3 & 2.25 & 2.03 \\
50\% & 16 & 1.14 & 2.00 & 3 & 3.51 & 2.12 \\
80\% & 17 & 1.4\hphantom{0} & 2.11 & 3 & 5.33 & 2.22 \\
90\% & 18 & 1.63 & 2.17 & 3 & 6.62 & 2.26 \\
\hline
\end{tabular*}
\end{table}

We have investigated the effects of sample size and magnitude of the
coefficients on PI.
The results show that (i) given the regression function and the noise level,
the value of PI indicates whether the problem is ``practically''
parametric/nonparametric at the current sample size; (2) given the noise
level and the sample size, when the nonparametric part is very weak,
PI has a large value, which properly indicates that the nonparametric
part is negligible; but as the
nonparametric part gets strong enough, PI will drop close to 1,
indicating a clear nonparametric scenario.
For a parametric scenario, the stronger the signal, the larger PI as is
expected.
See \cite{Liu2011} for details.
\end{Example}


\subsection{Multiple predictors}
In the multiple-predictor examples, we are going to do all subset
selection. We generate data from a linear model (except Example~\ref{Example7}):
$Y=\beta^{T}{\mathbf x}+\sigma\varepsilon$, where ${\mathbf x}$ is generated
from a multivariate normal distribution with mean 0, variance 1, and
correlation structure given in each example. For each generated data
set, we
apply the Branch and Bound algorithm \cite{Hand81} to do all subset
selection by BIC and then calculate the PI value (part of our code is
modified from the aster package of Geyer \cite{Geyer08}).
Unless otherwise stated, in these examples, the sample size is 200 and
we replicate 300 times.
The first two examples were used in \cite{Tibshirani96}.
\begin{Example}\label{Example3}
$\!\!\!\beta\,{=}\,(3,1.5,0,0,2,0,0,0)^{T}$. The correlation between $x_{i}$
and~$x_{j}$ is $\rho^{|i-j|}$ with $\rho=0.5$. We set $\sigma=5$.
\end{Example}
\begin{Example}\label{Example4}
Differences from Example \ref{Example3}: $\beta_{j}=0.85, \forall j$ and $\sigma=3$.
\end{Example}
\begin{Example}\label{Example5}
$\beta=(0.9, 0.9, 0, 0, 2, 0, 0, 1.6, 2.2, 0,0,0,0)^{T}$. There are 13
predictors and the correlation between $x_{i}$ and $x_{j}$ is $\rho
=0.6$ and $\sigma=3$.
\end{Example}
\begin{Example}\label{Example6}
This example is the same as Example \ref{Example5} except that $\beta=(0.85$, $0.85,
0, 0, 2, 0, 0, 0, 0, 0,0,0,0)^{T}$ and $\rho=0.5$.
\end{Example}
\begin{Example}\label{Example7}
This example is the same as Example \ref{Example3} except that we add a nonlinear
component in the mean function and $\sigma=3$, that is, $Y=\beta
^{T}{\mathbf x}+\phi(u) + \sigma\varepsilon$, where $u \sim \operatorname{uniform}(-4,4)$
and $\phi(u)=3(1-0.5u+2u^{2})e^{-u^{2}/4}$.
All subset selection is carried out with predictors $x_{1},\ldots,
x_{8}, u, \ldots, u^{8}$ which are coded as 1--8 and A--G in Table \ref
{multiple-predictor1}.

%
\begin{table}
\tablewidth=165pt
\caption{Proportion of selecting true model}
\label{multiple-predictor1}
\begin{tabular*}{\tablewidth}{@{\extracolsep{\fill}}llc@{}}
\hline
\textbf{Example} & \textbf{True model} &  \textbf{Proportion} \\
\hline
3 & 125  & 0.82 \\
4 & 12345678 & 0.12 \\
5 & 12589 & 0.43 \\
6 & 125 &  0.51 \\
7 & 1259ABCEG*   & 0.21 \\
\hline
\end{tabular*}
\vspace*{-6pt}
\end{table}

\begin{table}[t]
\tablewidth=150pt
\caption{Quartiles of PIs}
\label{multiple-predictor2}
\begin{tabular*}{\tablewidth}{@{\extracolsep{\fill}}lccc@{}}
\hline
\textbf{Example}   &  \textbf{Q1}    &  \textbf{Q2}   &  \textbf{Q3}      \\
\hline
3         &  1.26  &  1.51 &  1.81    \\
4         &  1.02  &  1.05 &  1.10    \\
5         &  1.05  &  1.15 &  1.35     \\
6         &  1.09  &  1.23 &  1.56     \\
7         &  1.02  &  1.07 &  1.16     \\
\hline
\end{tabular*}
\end{table}

The selection behaviors and PI values are reported in Tables \ref
{multiple-predictor1} and
\ref{multiple-predictor2}, respectively. From those results, we see
that the PIs are large for Example~\ref{Example3} and small for Example~\ref{Example4}. Note that
in Example~\ref{Example3} we have 82\% chance selecting the true model, while in
Example \ref{Example4} the chance is only 12\%. Although both Examples \ref{Example3} and \ref{Example4} are
of parametric nature, we would call Example \ref{Example4} ``practically
nonparametric'' in the sense that at the given sample size many models
are equally likely and the issue is to balance the approximation error
and estimation error. For Examples \ref{Example5} and \ref{Example6}, the PI values are
in-between, so are the chances of selecting the true models.
Note that the median PI values in Examples \ref{Example5} and \ref{Example6} are around 1.2.
These examples together show that the values of PI provide sensible
information on how strong the parametric message is and that
information is consistent with stability in selection.
\end{Example}

Example \ref{Example7} is quite interesting. Previously, without the $\phi(u)$
component, even at $\sigma=5$,
large values of PI are seen. Now with the nonparametric component
present, the PI values are close to 1. [The asterisk (*) in Table
\ref{multiple-predictor1} indicates the model is the most frequently
selected one instead of being the true model.]

More simulation results are given in \cite{Liu2011}. First, an
illuminating example shows that with
specially chosen coefficients, PI switches positions several times, as
they should, in declaring practical parametricness or nonparametricness
as more and more information is available.
Second, it is shown that PI is informative on reliability of inference
after model selection. When PI is large (Example \ref{Example3}), confidence
intervals based on the selected model are quite
trustworthy, but when PI is small (Example \ref{Example4}), the actual coverage
probability intended at 95\% is typically
around 65\%.
While it is now well known that model selection has an impact on subsequent
statistical inferences (see, e.g., \cite
{Zhang1990,HurvichTsai1990,Faraway1992,KabailaLeeb06}),
the value of PI can provide valuable information on the parametricness
of the underlying regression function and hence on how confident we are
on the accuracy of subsequent inferences.
Third, it is shown that an
adaptive choice between AIC and BIC based on the PI value (choose BIC
when PI is larger than 1.2) indeed
leads to nearly the better performance of AIC and BIC and thus beats
both AIC and BIC in an overall sense.
So PI provides helpful information regarding whether AIC or BIC works better
(or they have similar performances) in risks of estimation. Therefore,
PI can be viewed as a Performance Indicator of AIC versus BIC.

Based on our numerical investigations, in nested model problems (like
order selection for series
expansion), a cutoff point of $c=1.6$ seems proper. In subset
selection problems, since
the infimum in computing PI is taken over many models, the cutoff
point is expected to
be smaller, and $1.2$ seems to be quite good.

\section{Real data examples}\label{sec6}

In this section, we study three data sets: the Ozone data with 10
predictors and $n=330$ (e.g.,
\cite{Breiman85}),
the Boston housing data with 13 predictors and $n=506$ (e.g.,
\cite{Harrison78}), and the Diabetes data
with 10 predictors and $n=442$ (e.g., \cite{Efron04}).

In these examples, we conduct all subset selection by BIC using the
Branch and Bound algorithm. Besides finding the PI values for the full
data, we also do the same with sub-samples from the original data at
different sample sizes. In addition, we carry out a parametric
bootstrap from the model selected by BIC based on the original data to
assess the stability of model selection.

Based on sub-sampling at the sample size 400, we found that the PIs for
the ozone data are mostly larger
than 1.2, while those for the Boston housing data are smaller
than 1.2. Moreover, the parametric bootstrap suggests that for the
Ozone data, the model selected from the full data still reasonably
stands out even when the sample size is reduced to about 200 and noises are
added.
Similar to the simulation results in Section \ref{sec5}, by parametric
bootstrap at the original sample size from the selected model,
combining AIC and BIC based on PI shows good overall
performance in estimating the regression function. The combined
procedure has a statistical risk close to the better one of AIC and BIC
in each case. Details can be found in
\cite{Liu2011}.

\section{Conclusions}\label{sec7}

Parametric models have been commonly used to estimate a
finite-dimensional or infinite-dimensional
function. While there have been\vadjust{\goodbreak} serious debates on which model selection
criterion to use to choose a~candidate model
and there has been some work on combining the strengths of very
distinct model
selection methods, there is a major lack of understanding on statistically
distinguishing between
scenarios that favor one method (say AIC) and those that favor another
(say BIC).
To address this issue, we have derived a parametricness index (PI)
that has the desired theoretical property: PI
converges in probability to infinity for parametric scenarios and to~1
for nonparametric ones.
The use of a potentially consistent model
selection rule (i.e., it will be consistent if one of the candidate
models is true)
in constructing PI effectively prevents overfitting when we are
in a parametric scenario. The comparison of the selected model with a
subset model
separates parametric and nonparametric scenarios through the distinct
behaviors of
the approximation errors of these models in the two different situations.

One interesting consequence of the property of PI is that a choice
between AIC and
BIC based on its value ensures that the resulting regression estimator
of $f$ is automatically
asymptotically efficient for both parametric and nonparametric
scenarios, which
clearly cannot be achieved by any deterministic choice of the penalty
parameter $\lambda_n$ in the
criteria of the form $-\log\mbox{-}\mathrm{likelihood}+\lambda_n m_{k}$, where
$m_k$ is the number of
parameters in the model~$k$. Thus, an adaptive regression estimation to
simultaneously suit
parametric and nonparametric scenarios is realized through the
information provided by PI.

When working with parametric candidate models,
we advocate a practical view on parametricness/nonparametricness. In
our view,
a parametric scenario is one where a relatively parsimonious model reasonably
stands out. Otherwise, the selected model is most likely a tentative compromise
between goodness of fit and model complexity, and the recommended model
is most likely to change when the
sample size is slightly increased.

Our numerical results seem to be very encouraging. PI is informative,
giving the
statistical user an idea on how much one can trust the selected model
as the ``true''
one. When PI does not support the selected model as the ``right''
parametric model for the
data, we have demonstrated that estimation standard errors
reported from the selected model are often too small compared to the
real ones, that
the coverage of the resulting confidence intervals are much smaller
than the nominal levels,
and that mode selection uncertainty is high. In
contrast, when PI strongly endorses the selected model, model selection
uncertainty is much
less a concern and the resulting estimates and interpretation are
trustworthy to a large extent.

Identifying a stable and strong message in data as is expressed by a meaningful
parametric model, if existing,
is obviously important.
In biological and social sciences, especially observational studies, a
strikingly reliable parametric
model is often too much to ask for. Thus, to us, separating scenarios
where one model is reasonably
standing out and is expected to shine over other\vadjust{\goodbreak} models for sample
sizes not too much larger than
the current one from those where the selected model is simply the lucky
one to be chosen among multiple
equally performing candidates is an important step beyond simply
choosing a model based on one's
favorite selection rule or, in the opposite direction, not trusting any
post model selection
interpretation due to existence of model selection uncertainty.

For the other goal of regression function estimation, in application,
one typically applies a model
selection method, or considers estimates from two (or more) model
selection methods to see if they agree with each other. In light of PI
(or similar model selection diagnostic measures), the situation can be
much improved: one adaptively applies the better model selection
criterion to improve performance in estimating the regression function.
We have
focused on the competition between AIC and BIC, but similar measures
may be constructed for comparing other model selection methods that
are derived from different principles or under different assumptions.
For instance, the focused information criterion (FIC)
\cite{Claeskens03,Claeskens08} emphasizes
performance at a~given estimand, and it seems interesting to understand
when FIC improves over AIC and how to take advantages of both in an
implementable
fashion.\looseness=-1

For the purpose of estimating the regression function,
it has been suggested that AIC performs better for a nonparametric
scenario and BIC better for a parametric one (see \cite{Yang072} for a
study on the issue in a simple setting). This is asymptotically
justified but certainly not
quite true in reality. Our numerical results have demonstrated that for some
parametric regression functions, AIC is much better. On the other hand,
for an infinite-dimensional regression function, BIC can give a much
more accurate estimate. Our numerical results tend to suggest that
when PI is high and thus we are in a~practical
parametric scenario (whether the true regression function is
finite-dimensional or not), BIC tends to be better for regression
estimation; when PI is close to 1 and thus we are in a practical
nonparametric scenario, AIC tends to be better.

Finally, we point out some limitations of our work. First, our results
address only
linear models under Gaussian errors. Second, more understanding on the
choices of $\lambda_{n}$, $d$, and the best cutoff
value $c$ for PI is needed. Although the choices recommended in this
paper worked
very well for the numerical examples we have studied, different values
may be proper
for other situations (e.g., when the predictors are highly correlated and/or
the number of predictors is comparable to the sample size).


\begin{appendix}\label{app}
\section*{Appendix}

The following fact will be used in our proofs (see \cite{Yang99}).

\textit{Fact.} If $Z_{m}\sim\chi_{m}^{2}$, then
\begin{eqnarray*}
P(Z_{m}-m \ge\kappa m)&\le& e^{-{m}(\kappa-\ln(1+\kappa))/{2}}\qquad
\forall \kappa>0, \\
P(Z_{m}-m \le-\kappa m)&\le& e^{-{m}(-\kappa-\ln(1-\kappa
))/{2}}\qquad \forall 0<\kappa<1.\vadjust{\goodbreak}
\end{eqnarray*}

For ease of notation, we denote $P_{k^{(s)},k}=M_{k}-M_{k^{(s)}}$
by $P$, $\mathit{rem}_{1}(k)=e_{n}^{T}(f_{n}-M_{k}f_{n}) $ and
$\mathit{rem}_{2}(k)=\| (I_{n}-M_{k})e_{n}\|^{2}/\sigma^{2}-n $ in the
proofs. Then
%
%
\begin{eqnarray}
\label{eq8}
\bigl\|\bigl(I_{n}-M_{k^{(s)}}\bigr) e_{n}\bigr\|^{2}&=& \|(I_{n}-M_{k}) e_{n}\|^{2} +\|P
e_{n}\|^{2},
\\
\label{eq9}
\bigl\|\bigl(I_{n}-M_{k^{(s)}}\bigr) f_{n}\bigr\|^{2}&=& \|(I_{n}-M_{k}) f_{n}\|^{2} +\|P
f_{n}\|^{2},
\\
\label{eq10}
\mathit{rem}_{1}\bigl(k^{(s)}\bigr)&=& \mathit{rem}_{1}(k)+e_{n}^{T}P
f_{n}.
\end{eqnarray}

For the proofs of the theorems in the case of $\sigma$ known, without
loss of generality, we assume $\sigma^{2}=1$. In all the proofs, we
denote $IC_{\lambda_{n}, d}(k)$ by $ IC(k)$.\looseness=-1
\begin{pf*}{Proof of Theorem \ref{thm1} \textup{(parametric, $\sigma$
known)}}
Under the assumption that $P(\hat{k}_{n}=k_{n}^{*})\rightarrow1$, we
have $\forall\varepsilon>0$, $\exists n_{1}$ such that
$P(\hat{k}_{n}=k_{n}^{*})>1-\varepsilon$ for $n>n_{1}$.

Since
$ \|{\mathbf Y}_{n}-\hat{{\mathbf Y}}_{k}\|^{2}=\|(I_{n}-M_{k})f_{n}\|^{2}+\|
(I_{n}-M_{k})e_{n}\|^{2}
+2\,\mathit{rem}_{1}(k)$, for any $k_{n}^{*(s)}$ being a sub-model of $k_{n}^{*}$
with $r_{k_{n}^{*(s)}}=r_{k_{n}^{*}}-1 $,
we know that $\frac{IC(k_{n}^{*(s)})}{IC(k_{n}^{*})}$ is equal to
\begin{eqnarray*}
&&
\frac{\|{\mathbf Y}_{n}-\hat{{\mathbf Y}}_{k_{n}^{*(s)}}\|^{2}+\lambda_{n}\log
(n)r_{k_{n}^{*(s)}}-n+dn^{1/2}\log(n)}{\|{\mathbf Y}_{n}-\hat{{\mathbf Y}}_{k_{n}^{*}}\|^{2}+\lambda_{n}\log(n)r_{k_{n}^{*}}-n+dn^{1/2}\log
(n)} \\
&&\qquad=\bigl(\bigl\|\bigl(I_{n}-M_{k_{n}^{*(s)}}\bigr)f_{n}\bigr\|^{2}
+\mathit{rem}_{2}\bigl(k_{n}^{*(s)}\bigr)
+2\,\mathit{rem}_{1}\bigl(k_{n}^{*(s)}\bigr)\\
&&\hspace*{91.3pt}{}
+\lambda_{n}\log(n)(r_{k_{n}^{*}}-1)+dn^{
{1/2}}\log(n)\bigr)\\
&&\qquad\quad\hspace*{0pt}{}\times\bigl({ \mathit{rem}_{2}(k_{n}^{*})
+\lambda_{n}\log(n)r_{k_{n}^{*}} +dn^{1/2}\log(n)}\bigr)^{-1} .
\end{eqnarray*}
By the fact on $\chi^2$ distribution,
\begin{eqnarray*}
&& P\bigl(\|(I_{n}-M_{k_{n}^{*}})e_{n}\|^{2}-(n-r_{k_{n}^{*}}) \ge\kappa
(n-r_{k_{n}^{*}})\bigr)\\
&&\qquad\le e^{-({n-r_{k_{n}^{*}}})(\kappa-\ln(1+\kappa
))/{2}} \qquad\mbox{for } \kappa>0,
\\
&& P\bigl(\|(I_{n}-M_{k_{n}^{*}})e_{n}\|^{2}-(n-r_{k_{n}^{*}})\le-\kappa
(n-r_{k_{n}^{*}})\bigr)\\
&&\qquad\le e^{-({n-r_{k_{n}^{*}}})(-\kappa-\ln(1-\kappa
))/{2}} \qquad\mbox{for } 0<\kappa<1.
\end{eqnarray*}
For the given $\tau>0$, let $\kappa=\frac{n^{{1}/{2}+\tau
}h_{n}}{n-r_{k_{n}^{*}}}$ for some $h_{n} \rightarrow0$. Note that
when $n$ is large enough, say $n>n_{2}>n_{1}$, we have $0<\kappa=\frac
{n^{{1}/{2}+\tau}h_{n}}{n-r_{k_{n}^{*}}}<1$.
Since $
x-\log(1+x)\ge\frac{1}{4}x^{2}$ and $
-x-\log(1-x)\ge\frac{1}{4}x^{2} $ for $0<x<1$, we have
\begin{eqnarray*}
P\bigl(\bigl| \|(I_{n}-M_{k_{n}^{*}})e_{n}\|^{2}-(n-r_{k_{n}^{*}})
\bigr| \ge h_{n}n^{{1/2}+\tau} \bigr)
&\le& 2e^{-({n-r_{k_{n}^{*}}})\kappa^{2}/{8}} \le2 e^{-
n^{2\tau}h_{n}^{2}/{8}}.
\end{eqnarray*}

Since for $Z\sim N(0, 1)$, $\forall t>0$, $P(|Z|\ge t)\le
e^{-t^{2}/2}$, we know that $\forall c>0$,
\[
P\biggl(\frac{|\mathit{rem}_{1}(k_{n}^{*(s)})|}{\|(I_{n}-M_{k_{n}^{*(s)}})f_{n}\|
^{2}} \ge c \biggr)
\le e^{-c^{2}\|(I-M_{k_{n}^{*(s)}})f_{n}\|^{2}/2}.
\]
Thus, $ | \frac{IC(k_{n}^{*(s)})}{IC( k_{n}^{*})} |$ is no
smaller than
\begin{eqnarray*}
&&\bigl(\bigl|\bigl\|\bigl(I_{n}-M_{k_{n}^{*(s)}}\bigr)f_{n}\bigr\|^{2}
+\mathit{rem}_{2}\bigl(k_{n}^{*(s)}\bigr)\\
&&\qquad\hspace*{0pt}{}+2\,\mathit{rem}_{1}\bigl(k_{n}^{*(s)}\bigr)+\lambda_{n}\log(n)
(r_{k_{n}^{*}}-1)+dn^{
{1/2}}\log(n) \bigr|\bigr)\\
&&\qquad\hspace*{11pt}{}\times\bigl(h_{n}n^{1/2+\tau}+r_{k_{n}^{*}}\bigl(\lambda_{n}\log
(n)-1\bigr)+dn^{1/2}\log(n)\bigr)^{-1}
\end{eqnarray*}
%
with probability higher than $1-
2e^{-n^{2\tau}h_{n}^{2}/{8}}$.\vspace*{1pt}

Note that $IC(k_{n}^{*(s)}) $ is no smaller than
\[
(1-2c)\bigl\|\bigl(I_{n}-M_{k_{n}^{*(s)}}\bigr)f_{n}\bigr\|^{2}-h_{n}n^{1/2+\tau
}+(r_{k_{n}^{*}}-1)\bigl(\lambda_{n}\log(n)-1\bigr)+dn^{1/2}\log(n)
\]
with probability higher than $1-e^{-n^{2\tau
}h_{n}^{2}/{8}}-e^{-c^{2}\|(I-M_{k_{n}^{*(s)}})f_{n}\|^{2}/2}$.
Since $A_{n}$ is of order higher than $h_{n}n^{{1/2}+\tau} $ and
for $c<1/2$ (to be chosen), there exists $n_{3}>n_{2}$ such that
$IC(k_{n}^{*(s)}) $ is positive for $n>n_{3}$ and $|
\frac{IC(k_{n}^{*(s)})}{IC( k_{n}^{*})} |$ is no smaller than
\begin{eqnarray*}
&&\bigl((1-2c)\bigl\|\bigl(I_{n}-M_{k_{n}^{*(s)}}\bigr)f_{n}\bigr\|^{2}-h_{n}n^{1/2+\tau
}\\
&&\qquad\hspace*{0pt}{}+(r_{k_{n}^{*}}-1)\bigl(\lambda_{n}\log(n)-1\bigr)+dn^{1/2}\log(n)\bigr)
\\
&&\qquad\hspace*{11pt}{}\times\bigl(h_{n}n^{1/2+\tau}+r_{k_{n}^{*}}\lambda_{n}\log(n)+dn^{1/2}\log(n)\bigr)^{-1}
\end{eqnarray*}
%
with probability higher than
$1\,{-}\,2e^{-n^{2\tau}h_{n}^{2}/{8}}\,{-}\,(e^{-n^{2\tau}h_{n}^{2}/{8}}
\,{+}\,e^{-c^{2}\|(I-M_{k_{n}^{*(s)}})f_{n}\|^{2}/2})$.

Then for $n>n_{3}$, $\inf_{k^{*(s)}_{n}} |
\frac{IC(k_{n}^{*(s)})}{IC( k_{n}^{*})} | $ is lower bounded by
\[
\frac{(1-2c)A_{n}-h_{n}n^{1/2+\tau}+(r_{k_{n}^{*}}-1)(\lambda_{n}\log
(n)-1)+dn^{1/2}\log(n)}
{h_{n}n^{1/2+\tau}+r_{k_{n}^{*}}\lambda_{n}\log(n)+dn^{1/2}\log(n)}
\]
with probability higher than $ 1- 2e^{-n^{2\tau
}h_{n}^{2}/{8}}-r_{k_{n}^{*}}\cdot( e^{-n^{2\tau
}h_{n}^{2}/{8}}+e^{-c^{2}A_{n}/2})$.

According to conditions (P1) and (P2), $r_{k_{n}^{*}}=o(n^{
{1}/{2}+\tau})/(\lambda_{n}\log(n))$ and~$ A_{n} $ is of order
$n^{1/2+\tau}$ or higher, we can choose
$h_{n}$ such that $2e^{-n^{2\tau}h_{n}^{2}/{8}}
+r_{k_{n}^{*}}\cdot( e^{-n^{2\tau
}h_{n}^{2}/{8}}+e^{-c^{2}A_{n}/2}) \rightarrow0$.

For example, taking $h_{n}=n^{-\tau/3}$, then
\begin{eqnarray*}
\inf_{k^{*(s)}_{n}} \biggl|
\frac{IC(k_{n}^{*(s)})}{IC( k_{n}^{*})} \biggr| &\ge& \frac
{(1-2c)A_{n}-n^{1/2+2\tau/3}+(r_{k_{n}^{*}}-1)\lambda_{n}\log
(n)+dn^{1/2}\log(n)}
{n^{1/2+2\tau/3}+r_{k_{n}^{*}}\lambda_{n}\log(n)+dn^{1/2}\log(n)} \\
:\!& = & \mathrm{bound}_{n}
\end{eqnarray*}
with probability higher than $ 1-2e^{-n^{4\tau/3}/8}
-r_{k_{n}^{*}} (e^{-n^{4\tau/3}/8}+e^{-c^{2}A_{n}/2}) := 1-q_{n}$.

With $c<1/2$, $A_{n}$ of order $n^{1/2+\tau}$ or higher, and
$r_{k_{n}^{*}}\lambda_{n}\log(n)=o(A_{n})$, we have
that $\forall M >0, \exists n_{4}>n_{3}$ such that $\mathrm{bound}_{n}\ge M$
and $q_{n}\le\varepsilon$ for $n>n_{4}$.
Thus $\mathrm{PI}_{n}\stackrel{p}{\longrightarrow} \infty$.
\end{pf*}
\begin{pf*}{Proof of Theorem \ref{thm2} \textup{(nonparametric, $\sigma$
known)}}
Similar to the proof of Theorem \ref{thm1}, consider
$\frac{IC(\hat{k}_{n}^{(s)})}{IC(\hat{k}_{n})}$
for any $\hat{k}_{n}^{(s)}$ being a sub-model of $\hat{k}_{n}$ with one
fewer term, and we have
\begin{eqnarray*}
\frac{IC(\hat{k}_{n}^{(s)})}{IC( \hat{k}_{n})}&=&1+ \bigl({\|Pf_{n}\|
^{2}+\|Pe_{n}\|^{2}
+e_{n}^{T}Pf_{n}-\lambda_{n}\log(n)}\bigr)\\
&&\hspace*{16pt}{}\times\bigl(\|(I_{n}-M_{\hat{k}_{n}})f_{n}\|^{2}
+\mathit{rem}_{2}(\hat{k}_{n})\\
&&\qquad\hspace*{10pt}{}+2\,\mathit{rem}_{1}(\hat{k}_{n})
+\lambda_{n}\log(n)r_{\hat{k}_{n}}+dn^{1/2}\log(n)\bigr)^{-1}.
\end{eqnarray*}
Next, consider the terms in the above equation for any model $k_{n}$.
For ease of notation, we write
$ B_{r_{k_{n}},n}=B_{r_{k_{n}}}$, where $r_{k_{n}}$ is the rank of the
projection matrix of model $k_{n}$.

As in the proof of Theorem \ref{thm1},
$\forall c_{1}>0$,
\begin{eqnarray*} \label{equation6}
&& P\biggl(\frac{|\mathit{rem}_{1}(k_{n})|}{(\lambda_{n}\log(n)-1)r_{k_{n}}+\|
(I_{n}-M_{k_{n}})f_{n}\|^{2}+dn^{1/2}\log(n)} \ge c_{1} \biggr)
\nonumber\\
&&\qquad\le
e^{-c_{1}^{2}({(\lambda_{n}\log(n)-1)r_{k_{n}}+\|
(I_{n}-M_{k_{n}})f_{n}\|^{2}+dn^{1/2}\log(n)})/{2}} \le
e^{-c_{1}^{2}B_{r_{k_{n}}}/2} .
\end{eqnarray*}

Similarly, $\forall c_{2}>0$,
%
%
\begin{eqnarray}\label{equation7}\qquad
P\biggl(\frac{|e_{n}^{T}Pf_{n}|}{B_{r_{k_{n}}}}\ge c_{2}\biggr) &\le&
e^{-{c_{2}^{2}B_{r_{k_{n}}}^{2}}/({2\|Pf_{n}\|^{2}})}\nonumber\\[-8pt]\\[-8pt]
&\le&
e^{-c_{2}^{2}B_{r_{k_{n}}}/2} \qquad(\mbox{if } \|Pf_{n}\|
^{2}\le B_{r_{k_{n}}}),\nonumber
\\
%
%
\label{equation8}
P\biggl(\frac{|e_{n}^{T}Pf_{n}|}{\|Pf_{n}\|^{2}}\ge c_{2}\biggr)
&\le& e^{-{c_{2}^{2}\|Pf_{n}\|^{2}}/{2}} \nonumber\\[-8pt]\\[-8pt]
&\le&
e^{-c_{2}^{2}B_{r_{k_{n}}}/2} \qquad(\mbox{if } \|Pf_{n}\|^{2}>
B_{r_{k_{n}}}).\nonumber
\end{eqnarray}
%

Also,
\[
P\bigl(\|(I_{n}-M_{k_{n}})e_{n}\|^{2}-(n-r_{k_{n}})\le-\kappa
(n-r_{k_{n}})\bigr)\le
e^{-({n-r_{k_{n}}})(-\kappa-\log(1-\kappa))/{2}} .
\]
We can choose $\kappa$ such that $\kappa
(n-r_{k_{n}})=\gamma B_{r_{k_{n}}}$ for some
$0< \gamma<1 $. Note that $-x-\log(1-x)> x^{2}/2$ for $0<x<1$.
Then
%
%
\begin{equation} \label{equation9}
P\bigl(\|(I_{n}-M_{k_{n}})e_{n}\|^{2}-(n-r_{k_{n}})\le-\gamma_{n}
B_{r_{k_{n}}}\bigr)\le
e^{-{\gamma^{2} B_{r_{k_{n}}}^{2}}/({4(n-r_{k_{n}})})} .
\end{equation}
For a sequence $D_{n}>0$ (to be chosen), we have
\[
P(\|Pe_{n}\|^{2}-1\ge D_{n})\le e^{-(D_{n} -\log(1+D_{n}
))} .
\]
For\vadjust{\goodbreak} $x>1$, $x-\log(1+x)>x/2$. So $ P(\|Pe_{n}\|^{2}-1\ge
D_{n})\le e^{-D_{n}/2 }$ for $D_{n}>1$.

Since $\hat{k}_{n}$ is random, we apply union bounds on the exception
probabilities. According to condition (N1), for\vspace*{1pt} any $\varepsilon>0$,
there exists $ n_{1} $ such that $P(a_{n}\le r_{\hat{k}_{n}}\le
b_{n})\ge1-\varepsilon$ for $n>n_{1}$.
As will be seen, when $n$ is large enough, the following quantities can
be arbitrarily small for
appropriate choice of $\gamma$, $D_{n}$, $c_{1}$
and $c_{2}$:
\begin{eqnarray*}
&&\sum_{j=a_{n}}^{b_{n}} N_{j} \cdot e^{-{\gamma^{2} B_{j,n}^{2}
}/({4(n-j)})},\qquad
\sum_{j=a_{n}}^{b_{n}} N_{j}\cdot L_{j} \cdot e^{-D_{n}/2 } ,\\[-3pt]
&&\sum_{j=a_{n}}^{b_{n}} N_{j} \cdot e^{-c_{1}^{2}B_{j,n}/2},\qquad
\sum_{j=a_{n}}^{b_{n}} N_{j}\cdot L_{j} \cdot
e^{-c_{2}^{2}B_{j,n}/2} .
\end{eqnarray*}
More precisely, we claim that there exists $ n_{2}>n_{1} $ such that
for $
n\ge n_{2}$,
%
%
\begin{eqnarray}\label{series}
&&\sum_{j=a_{n}}^{b_{n}} \bigl\{N_{j} \cdot\bigl(e^{-{\gamma^{2}
B_{j,n}^{2}}/({4(n-j)})} +
e^{-c_{1}^{2}B_{j,n}/2}\bigr)\nonumber\\[-10pt]\\[-10pt]
&&\qquad\hspace*{9.4pt}{} +
N_{j}\cdot L_{j} \cdot( e^{-D_{n}/2 } +
e^{-c_{2}^{2}B_{j,n}/2} )\bigr\} \le
\varepsilon.\nonumber
\end{eqnarray}

Then for $n>n_{2}$ with probability higher than $1-2\varepsilon$,
\begin{eqnarray*}
&a_{n}\le r_{\hat{k}_{n}} \le b_{n},&
\\[-2pt]
&\|(I_{n}-M_{\hat{k}_{n}})e_{n}\|^{2}-(n-r_{\hat{k}_{n}}) \ge-\gamma
B_{r_{\hat{k}_{n}}},&
\\[-2pt]
&\bigl\|P_{\hat{k}_{n}^{(s)},\hat{k}_{n}}e_{n}\bigr\|^{2} \le1+D_{n},&
\\[-2pt]
&|\mathit{rem}_{1}(\hat{k}_{n})| \le c_{1}\bigl(\bigl(\lambda_{n}\log(n)-1\bigr)r_{\hat
{k}_{n}}+\|(I_{n}-M_{\hat{k}_{n}})f_{n}\|^{2}+dn^{1/2}\log(n)\bigr),&
\\[-2pt]
&\bigl|e_{n}^{T}P_{\hat{k}_{n}^{(s)},\hat{k}_{n}}f_{n}\bigr| \le c_{2} B_{r_{\hat
{k}_{n}}} \quad\mbox{or}\quad
\bigl|e_{n}^{T}P_{\hat{k}_{n}^{(s)},\hat{k}_{n}}f_{n}\bigr| \le
c_{2}\bigl\|P_{\hat{k}_{n}^{(s)},\hat{k}_{n}}f_{n}\bigr\|^{2} .&
\end{eqnarray*}
Note that
%
%
\begin{eqnarray} \label{PI}
\mathrm{PI}_{n}
&=& 1 + \inf_{\hat{k}_{n}^{(s)}} \bigl(\bigl({\|Pf_{n}\|^{2}+\|Pe_{n}\|^{2}
+e_{n}^{T}Pf_{n}-\lambda_{n}\log(n)}\bigr)\nonumber\\[-2pt]
&&\qquad\hspace*{11.6pt}{}\times
\bigl(\|(I_{n}-M_{\hat{k}_{n}})f_{n}\|^{2}+\mathit{rem}_{2}(\hat{k}_{n})\\[-2pt]
&&\qquad\hspace*{11.6pt}\hspace*{16.7pt}{}+2\,\mathit{rem}_{1}(\hat{k}_{n})
+\lambda_{n}\log(n)r_{\hat{k}_{n}}+dn^{1/2}\log(n)\bigr)^{-1}\bigr).\nonumber
\end{eqnarray}

Also with\vspace*{1pt} probability higher than $1- 2\varepsilon$, the denominator in
(\ref{PI}) is bigger than
$(1-2c_{1})[\|(I_{n}-M_{\hat{k}_{n}})f_{n}\|^{2}+(\lambda_{n}\log
(n)-1)r_{\hat{k}_{n}}+dn^{1/2}\log(n)]
-\gamma B_{r_{\hat{k}_{n}}}$. Thus, when $2c_{1}+\gamma<1$, the
denominator in (\ref{PI}) is positive.

Then for $n>n_{2}$, with probability at $1-2\varepsilon$ we have
\begin{eqnarray*}
\mathrm{PI}_{n}&=& 1 + \Bigl({\inf_{\hat{k}_{n}^{(s)}} \bigl(\|Pf_{n}\|^{2}+\|Pe_{n}\|^{2}
+e_{n}^{T}Pf_{n}-\lambda_{n}\log(n)\bigr)}\Bigr)\\[-2pt]
&&\qquad\hspace*{-6.5pt}{}\times\bigl(\|(I_{n}-M_{\hat{k}_{n}})f_{n}\|^{2}
+\mathit{rem}_{2}(\hat{k}_{n})\\[-2pt]
&&\qquad\hspace*{10.3pt}{}+2\,\mathit{rem}_{1}(\hat{k}_{n})
+\lambda_{n}\log(n)r_{\hat{k}_{n}}+dn^{1/2}\log(n)\bigr)^{-1}.
\end{eqnarray*}
For $n>n_{2}$ with probability higher than $1- 2\varepsilon$, if $ \|
Pf_{n}\|^{2} \le B_{r_{\hat{k}_{n}}} $, then
\[
\mathrm{PI}_{n} -1
\le \frac{\inf_{\hat{k}_{n}^{(s)}} \|Pf_{n}\|^{2}+ 1+D_{n} +c_{2}
B_{r_{\hat{k}_{n}}}+\lambda_{n}\log(n)}
{(1-2c_{1}-\gamma)((\lambda_{n}\log(n)-1)r_{\hat{k}_{n}}+\|
(I_{n}-M_{\hat{k}_{n}})f_{n}\|^{2}+dn^{1/2}\log(n)) }
\]
and
\[
\mathrm{PI}_{n}-1 \ge \frac{\inf_{\hat{k}_{n}^{(s)}} \|Pf_{n}\|
^{2}- 1-D_{n} -c_{2}
B_{r_{\hat{k}_{n}}}-\lambda_{n}\log(n)}
{(1-2c_{1}-\gamma)((\lambda_{n}\log(n)-1)r_{\hat{k}_{n}}+\|
(I_{n}-M_{\hat{k}_{n}})f_{n}\|^{2}+dn^{1/2}\log(n)) } ,
\]
otherwise,
\[
\mathrm{PI}_{n} -1 \le \frac{\inf_{\hat{k}_{n}^{(s)}}
\|Pf_{n}\|^{2}+ 1+D_{n} +c_{2}
\|Pf_{n}\|^{2}+\lambda_{n}\log(n)}
{(1-2c_{1}-\gamma)((\lambda_{n}\log(n)-1)r_{\hat{k}_{n}}+\|
(I_{n}-M_{\hat{k}_{n}})f_{n}\|^{2}+dn^{1/2}\log(n)) }
\]
and
\[
\mathrm{PI}_{n} -1 \ge \frac{\inf_{\hat{k}_{n}^{(s)}} \|
Pf_{n}\|^{2}- 1-D_{n} -c_{2}
\|Pf_{n}\|^{2}-\lambda_{n}\log(n)}
{(1-2c_{1}-\gamma)((\lambda_{n}\log(n)-1)r_{\hat{k}_{n}}+\|
(I_{n}-M_{\hat{k}_{n}})f_{n}\|^{2}+dn^{1/2}\log(n)) }.
\]

Next, we focus on the case $\|Pf_{n}\|^{2} \le B_{r_{\hat{k}_{n}}}$.
The case of $\|Pf_{n}\|^{2} > B_{r_{\hat{k}_{n}}} $ can be similarly
handled. Note that $\sup_{a_{n}\le j\le b_{n} }\frac
{B_{j,n}}{n-j}:=\zeta_{n}^{\prime} \rightarrow0 $. Let $\zeta_{n}^{\prime\prime}=\zeta
_{n}+\zeta_{n}^{\prime} $.
Taking $ \gamma=\sqrt{4/5}, D_{n}=4\zeta_{n}^{\prime\prime}B_{r_{k_{n}}},
c_{2}=2\sqrt{\zeta_{n}^{\prime\prime}}, 0<c_{1}<\frac{1-\gamma}{2} $, then
\begin{eqnarray*}
&& \mathrm{PI}_{n}-1 \\
&&\qquad \le \frac{\inf_{\hat{k}_{n}^{(s)}} \|Pf_{n}\|^{2}+
1+ 4\zeta_{n}^{\prime\prime}B_{r_{\hat{k}_{n}}} + 2\sqrt{\zeta_{n}^{\prime\prime}}B_{r_{\hat
{k}_{n}}} +\lambda_{n}\log(n)}
{(1-2c_{1}-\gamma)((\lambda_{n}\log(n)-1)r_{\hat{k}_{n}}+\|
(I_{n}-M_{\hat{k}_{n}})f_{n}\|^{2}+dn^{1/2}\log(n)) }\\
&&\qquad\le \sup_{ a_{n}\le r_{k_{n}}\le b_{n} }\Bigl(\Bigl(\inf_{k_{n}^{(s)}} \|
Pf_{n}\|^{2}+ 1+ 4\zeta_{n}^{\prime\prime}B_{r_{k_{n}}} + 2\sqrt{\zeta
_{n}^{\prime\prime}}B_{r_{k_{n}}} +\lambda_{n}\log(n)\Bigr)\\
&&\qquad\hspace*{53.7pt}{}\times
\bigl((1-2c_{1}-\gamma)\bigl(\bigl(\lambda_{n}\log(n)-1\bigr)r_{k_{n}}\\
&&\qquad\hspace*{135pt}{}+\|
(I_{n}-M_{k_{n}})f_{n}\|^{2}+dn^{1/2}\log(n)\bigr)\bigr)^{-1}\Bigr)\\
&&\qquad:=  \mathrm{Upperbound}_{n}\\
&&\qquad\rightarrow 0\qquad  \mbox{according to (N3) and the fact that }
\zeta_{n}^{\prime\prime}\rightarrow0 \mbox{ as } n \rightarrow\infty.
\end{eqnarray*}
Similarly,
\begin{eqnarray*}
&&\mathrm{PI}_{n}-1\\
&&\qquad \ge -\frac{ 1+ 4\zeta_{n}^{\prime\prime}B_{r_{k_{n}}} + 2\sqrt{\zeta
_{n}^{\prime\prime}}B_{r_{k_{n}}} +\lambda_{n}\log(n)}
{(1-2c_{1}-\gamma)((\lambda_{n}\log(n)-1)r_{\hat{k}_{n}}+\|
(I_{n}-M_{\hat{k}_{n}})f_{n}\|^{2}+dn^{1/2}\log(n))}\\
&&\qquad\ge - \sup_{a_{n}\le r_{k_{n}}\le b_{n} }\bigl(\bigl({ 1+ 4\zeta
_{n}^{\prime\prime}B_{r_{k_{n}}} + 2\sqrt{\zeta_{n}^{\prime\prime}}B_{r_{k_{n}}} +\lambda
_{n}\log(n)}\bigr)\\
&&\qquad\hspace*{62.3pt}{}\times\bigl((1-2c_{1}-\gamma)\bigl(\bigl(\lambda_{n}\log(n)-1\bigr)r_{k_{n}}\\
&&\qquad\hspace*{145pt}{}+\|
(I_{n}-M_{k_{n}})f_{n}\|^{2}+dn^{1/2}\log(n)\bigr) \bigr)^{-1}\bigr)\\
&&\qquad:=\mathrm{Lowerbound}_{n} \\
&&\qquad\rightarrow 0 \qquad\mbox{according to  (N3) and the fact
that } \zeta_{n}^{\prime\prime}\rightarrow0.
\end{eqnarray*}
Therefore, $\forall\delta>0, \exists n_{3}$ such
that $\mathrm{Upperbound}_{n}\le\delta$ and $\mathrm{Lowerbound}_{n}\ge-\delta$ for
$n>n_{3}$. Thus,
$\forall\varepsilon>0,\delta>0, \exists N=\max(n_{2},n_{3} )$ such that
$P(| \mathrm{PI}_{n}-1 |\le\delta)\ge1-2\varepsilon$ for $n>N$. That is,
$\mathrm{PI}_{n}\stackrel{p}{\rightarrow} 1. $

To complete the proof, we just need to check the claim of (\ref
{series}). By condition (N2), $\forall\varepsilon>0, \exists
n_{\varepsilon} $ such that for $n\ge n_{\varepsilon}$, $\sum
_{j=a_{n}}^{b_{n}}c_{0} \cdot e^{-{
B_{j,n}^{2}}/({10(n-j)})}<\varepsilon
/4 $. Then for $n>n_{\varepsilon}$,
\begin{eqnarray*}
\sum_{j=a_{n} }^{b_{n}} N_{j} \cdot e^{-({\gamma^{2}
B_{j,n}^{2}})/({4(n-j)})} &\le& \sum_{j=a_{n} }^{b_{n}} c_{0}\cdot e^{
{ B_{j,n}^{2}}/({10(n-j)})} \cdot e^{-{\gamma^{2}
B_{j,n}^{2}}/({4(n-j)})} \\
&\le&\sum_{j=a_{n}}^{b_{n}}c_{0} \cdot e^{-{
B_{j,n}^{2}}/({10(n-j)})}<\varepsilon/4, \\
\sum_{j=a_{n} }^{b_{n}} N_{j}\cdot L_{j} \cdot e^{-D_{n}/2 }
&=& \sum_{j=a_{n} }^{b_{n}} N_{j} \cdot L_{j}
\cdot e^{- 2\zeta_{n}^{\prime\prime} B_{j,n} } \\
&\le&\sum_{j=a_{n} }^{b_{n}} c_{0}
\cdot e^{- \zeta_{n}^{\prime\prime} B_{j,n} } <\frac{\varepsilon}{4}.
\end{eqnarray*}
Similarly,
\[
\sum_{j=a_{n} }^{b_{n}} N_{j} \cdot e^{-c_{1}^{2}B_{j,n}/2} <\frac
{\varepsilon}{4} ,\qquad \sum_{j=a_{n} }^{b_{n}} N_{j}\cdot L_{j}
\cdot e^{-c_{2}^{2}B_{j,n}/2} < \frac{\varepsilon}{4}.
\]
Thus, claim (\ref{series}) holds and this completes the proof.
\end{pf*}

The proofs of the cases with unknown $\sigma$ in Theorems \ref{thm1} and \ref{thm3} are
almost the same as those
when $\sigma$ is known. Due to space limitation, we omit the details.
\end{appendix}

\section*{Acknowledgments}
The authors thank Dennis Cook, Charles Geyer, Wei Pan, Hui Zou and the
participants at a seminar that one of the authors gave in Department of
Statistics at Yale University for helpful comments and discussions.
Comments from all the reviewers, Associate Editors and Editors are appreciated.

\begin{supplement}[id=suppA]
\stitle{Details and more numerical examples}
\slink[doi]{10.1214/11-AOS899SUPP} 
\sdatatype{.zip}
\sfilename{aos899\_supp.zip}
\sdescription{We provide
complete descriptions and more results of our numerical work.}
\end{supplement}

%

%
\printaddresses

\end{document}